\documentclass[12pt]{amsart}
\usepackage{graphics,color,epsfig}
\usepackage{amsfonts,amssymb}

\newcommand{\C}{\mathbb C}
\newcommand{\HH}{\mathbb H}
\newcommand{\N}{\mathbb N}

\newcommand{\R}{\mathbb R}
\newcommand{\Z}{\mathbb Z}

\newcommand{\bc}{\mathbf{c}}
\newcommand{\bx}{\mathbf{x}}

\newcommand{\cB}{\mathcal B}
\newcommand{\cC}{\mathcal C}
\newcommand{\cU}{\mathcal U}

\newcommand{\tu}{\tilde u}
\newcommand{\tv}{\tilde v}

\DeclareMathOperator{\diag}{diag}
\DeclareMathOperator{\diam}{diam}
\DeclareMathOperator{\Hdim}{dim_H}
\DeclareMathOperator{\SL}{SL}

\newcommand{\eps}{\varepsilon}
\newcommand{\vhi}{\varphi}
\newcommand{\pmat}[1]{\begin{pmatrix}#1\end{pmatrix}}
\renewcommand{\Re}{\textrm{Re~}}
\renewcommand{\Im}{\textrm{Im~}}

\newtheorem{theorem}{Theorem}[section]
\newtheorem{definition}[theorem]{Definition}
\newtheorem{lemma}[theorem]{Lemma}

\newtheorem{corollary}[theorem]{Corollary}

\title[Hausdorff dimension for divergence in a product space]
{Hausdorff dimension of the set of points on divergent trajectories 
       of a homogeneous flow on a product space}
\author{Yitwah Cheung}
\address{San Francisco State University \\ 
San Francisco, CA 94132, U.S.A.}
\email{cheung@math.sfsu.edu}
\date{\today}

\begin{document}
\begin{abstract}
In this paper we compute the Hausdorff dimension of the set $D(\vhi_n)$ 
  of points on divergent trajectories of the homogeneous flow $\vhi_n$ 
  induced by the one-parameter subgroup $\diag(e^t,e^{-t})$ acting by 
  left multiplication on the product space $G^n/\Gamma^n$, where 
  $G=\SL(2,\R)$ and $\Gamma=\SL(2,\Z)$.  
We prove that $\Hdim D(\vhi_n)=3n-\frac{1}{2}$ for $n\ge2$.  
\end{abstract}
\maketitle

\section{Introduction}
Let $G=\SL(2,\R)$ be the special linear group of two-by-two matrices 
  with real entries and determinant one and let $\Gamma=\SL(2,\Z)$ be 
  the discrete subgroup formed by those with integer entries.  
Let $G(n)=G_1\times\dots\times G_n$ where $G_i=G$, 
  $\Gamma(n)=\Gamma_1\times\dots\times\Gamma_n$ where $\Gamma_i=\Gamma$,  
  and consider the noncompact homogeneous space 
  $$G(n)/\Gamma(n)=(G_1/\Gamma_1)\times\dots\times (G_n/\Gamma_n)$$ 
  and the flow induced by the one-parameter subgroup $\vhi_n:\R\to G(n)$, 
  $$t\to(g_t,\dots,g_t) \quad\text{where}\quad 
                        g_t=\pmat{e^{t/2}&0\\0&e^{-t/2}},$$  
  acting by left multiplication.  
The forward trajectory $(x_t)_{t\ge0}$ of a point $x\in G(n)/\Gamma(n)$ 
  is said to be \emph{divergent} if it eventually leaves every compact set, 
  i.e. for any compact subset $K\subset G(n)/\Gamma(n)$ there is a time $T$ 
  such that $x_t\not\in K$ for all $t>T$.  
Let $$D(\vhi_n):=
  \{x\in G(n)/\Gamma(n): x_t\to\infty\quad\text{ as }\quad t\to\infty\}$$  
  be the set of points whose forward trajectories are divergent; we note that 
  $D(\vhi_n)$ is the union of all forward divergent trajectories of $\vhi_n$.  
The Hausdorff dimension of a subset of $G(n)/\Gamma(n)$ is defined with respect 
  to any metric induced by a right invariant metric on $G(n)$, the choice being 
  irrelevant since Hausdorff dimension depends only on the Lipschitz class of 
  a metric.  

In this paper we compute the Hausdorff dimension of the set $D(\vhi_n)$.  
\begin{theorem}\label{thm:main}
  For any $n\ge2$, $\Hdim D(\vhi_n)=\dim G(n)-\frac{1}{2}$.  
\end{theorem}

Theorem~\ref{thm:main} was motivated by certain analogies between partially 
  hyperbolic homogeneous flows on finite volume noncompact spaces and the 
  Teichm\"uller flow on moduli spaces of holomorphic quadratic differentials.  
The Teichm\"uller flow can be defined as the restriction to the diagonal 
  subgroup $(\phi_t)$ of a certain action of $\SL(2,\R)$ on moduli space.  
In this context, Masur \cite{Ma} showed that for any $\SL(2,\R)$-orbit $X$ 
  the Hausdorff dimension of the set of points on divergent trajectories 
  of the Teichm\"uller flow is at most $\frac{1}{2}$; moreover, for a 
  generic $\SL(2,\R)$-orbit, it was shown by Masur-Smillie \cite{MS} that 
  the Hausdorff dimension is in fact positive.  

For any partially hyperbolic homogeneous flow $\vhi$ on a finite volume 
  noncompact homogeneous space, Dani \cite{Da} showed that the set $D(\vhi)$ 
  contains a countable union of submanifolds of positive codimension that 
  consists of points lying on \emph{degenerate} divergent trajectories.  
Moreover, he showed that in the case of $\R$-rank one, all divergent 
  trajectories are degenerate, while in the higher rank situation there 
  are \emph{nondegenerate} divergent trajectories.  (See also \cite{We} 
  for related questions pertaining to the notion of degeneracy.)  
Theorem~\ref{thm:main} shows that in the case of $\vhi_n$ for $n\ge2$ the 
  set of points that lie on degenerate divergent trajectories form a subset 
  of $D(\vhi_n)$ of positive Hausdorff codimension, supporting the idea that 
  nondegenerate divergent trajectories are more abundant than degenerate ones.  
We also note that for $n=1$, we have $\Hdim D(\vhi_1)=\dim G-1=2$.  

For further results concerning the trajectories (bounded or divergent) of 
  partially hyperbolic homogeneous flows and their applications to number 
  theory we refer the reader to \cite{Bu}, \cite{Kl}, \cite{KM}, \cite{St}, 
  and \cite{We}.  

To obtain the upper bound in Theorem~\ref{thm:main} we shall in fact show 
  that for some sufficiently large compact set $K\subset G(n)/\Gamma(n)$ 
  the points whose forward trajectory under eventually stays outside $K$ 
  form a set $S=S(K)$ of positive Hausdorff codimension.  
The compact set will depend on a parameter $\delta>0$ such that the upper 
  bound on $\Hdim S(K)$ tends to $3n-\frac{1}{2}$ as $\delta\to0$.  

\textit{Outline.}
In \S\ref{S:Endpt}, we use standard methods to reduce the computation 
  of the Hausdorff dimension of $D(\vhi_n)$ to that of the set $E_n^*$ 
  of endpoints of nondegenerate divergent trajectories.  
The set $E_n^*$ may naturally be thought of as a subset of $\R^n$ and 
  in \S\ref{S:PL} we give a characterisation of this set in terms of 
  an encoding that uses continued fractions.  
Using this encoding, we compute the lower bound on $\Hdim D(\vhi_n)$ 
  in \S\ref{S:Lower}.  
Then we introduce the notion of a self-similar covering in \S\ref{S:Selfsim} 
  as a convenient device for presenting the upper bound calculation, 
  which is presented in \S\ref{S:Upper}.  

\textit{Acknowledgments.}
The author would like to thank the referee for many useful comments and 
  suggestions on an earlier version of this paper, which went under a 
  slightly different title.  
The author would also like to thank Alex Eskin and Barak Weiss for their 
  interest in this problem.  
Last, but not least, the author is deeply indebted to his wife Ying Xu 
  for her constant and unwavering support.

\section{Endpoints of divergent trajectories}\label{S:Endpt}
In this section we consider the set of endpoints of divergent trajectories 
  and identify it with a subset of $\R^n$.  
We assume the notation already established in the introduction.  

\textbf{Ideal boundary and Bruhat decomposition.}
The forward trajectories of two points in $G(n)$ are asymptotic if and 
  only if they belong to the same right coset of the parabolic subgroup 
  $P(n)=P_1\times\dots\times P_n$ where $P_i=P$ for $i=1,\dots,n$ and 
  $$P=\{p\in G: g_tpg_{-t} \text{ stays bounded as }t\to\infty \}.$$  
Thus, the set $D(\vhi_n)$ is a union of right cosets of $P(n)$.  

The ideal boundary, whose points are asymptotic classes of trajectories 
  of $\vhi_n$, is represented by the right coset space $$P(n)\backslash 
     G(n)=(P_1\backslash G_1)\times\dots\times(P_n\backslash G_n).$$  
The unipotent radicals of $P$ and $P(n)$ are respectively given by 
  $$N=\{u\in G: g_tug_{-t}\to e \text{ as }t\to\infty \}$$ 
  and $N(n)=N_1\times\dots\times N_n$ where $N_i=N$.  
The Bruhat decomposition establishes a one-to-one correspondence between 
  right cosets of $P(n)$ that differ from $P(n)$ itself with elements 
  of $N(n)$:  $$G(n)- P(n)=\bigcup_{u\in N(n)} P(n)wu$$ 
  where $w=(w_1,\dots,w_n)$ and $w_i=\pmat{0&-1\\1&0}$ for $i=1,\dots,n$.  
Let $E_n$ be the subset of $N(n)$ such that 
  $$(D(\vhi_n)- P(n))=\bigcup_{u\in E_n} P(n)wu.$$  

\textbf{Rational points.}
A right coset of $P_i$ is a $\Gamma_i$-\emph{rational point} if its 
  stabiliser under the action of $G_i$ by right multiplication contains 
  a maximal parabolic subgroup of $\Gamma_i$.  
Similarly, a right coset of $P(n)$ is $\Gamma(n)$-\emph{rational} if 
  its stabiliser under the $G(n)$ action on $P(n)\backslash G(n)$ 
  contains a maximal parabolic subgroup of $\Gamma(n)$.  
We note that $P(n)$ is a $\Gamma(n)$-rational point and also that the 
  set of $\Gamma(n)$-rational points forms a single orbit under the 
  action of $\Gamma(n)$ on $P(n)\backslash G(n)$ by right multiplication.  

We shall say a right coset of the form $P(n)wu$ is \emph{rational in 
  the $i$th coordinate} if its stabiliser under the action of $G(n)$ 
  on $P(n)\backslash G(n)$ contains the subgroup of $G(n)$ consisting 
  of elements of the form $(x_1,\dots,x_n)$ where $x_j$ is the identity 
  of $G_j$ for $j\neq i$ and $x_i\in P_i$.  
Note that $P(n)wu$ is $\Gamma(n)$-rational if and only if it is rational 
  in the $i$th coordinate for $i=1,\dots,n$.  
We say $P(n)wu$ is \emph{totally irrational} if it is not rational in 
  the $i$th coordinate for any $i\in\{1,\dots,n\}$.  
Let $E_n'$ be the set of points in $E_n$ corresponding to right cosets 
  of $P(n)$ that a totally irrational.  

We shall identify $N(n)$ with $\R^n$ and $E_n'$ with the corresponding 
  subset of $\R^n$.  
Since $\dim P(n)=2n$, the product formula for Hausdorff dimension (see 
  (\ref{eqn:product}) in \S\ref{S:Lower}) gives 
                $$\Hdim (P(n)E_n')=2n + \Hdim E_n'.$$ 
Since $D(\vhi_n)- P(n)E_n'$ is contained in a countable union 
  of submanifolds of $G(n)$ of dimension strictly less than $\dim G(n)$, 
  its Hausdorff dimension is bounded above by $3n-1$ 
  (see Lemma~\ref{lem:sup:Hdim} in \S\ref{S:Selfsim}.)  
Therefore, the proof of Theorem~\ref{thm:main} reduces to the statement 
  $$\Hdim E_n'=n-\frac{1}{2}\quad\text{for}\quad n\ge2.$$  

The identification of $N(n)$ with $\R^n$ can be made explicit as follows.  
The action of $G$ on the upper half plane $\HH^2=\{z\in\C: \Im z>0\}$ 
  by M\"obius transformations extends continuously to the boundary 
  $\R\cup\{\infty\}$; we assume this action is given as a right action.  
The boundary is naturally identified with $P\backslash G$ such that 
  the right coset $P$ corresponds to the ideal point $\infty$.  
Similarly, the group $G(n)$ acts on the product space $(\HH^2)^n$ and 
  the ideal boundary $P(n)\backslash G(n)$ is naturally identified with 
  the topological boundary of $(\HH^2)^n$ as a subset of $(\Hat\C)^n$ 
  where $\Hat\C=\C\cup\{\infty\}$.  
Under this identification, the set of right cosets of the form $P(n)wu$ 
  with $u\in N(n)$ corresponds to the subset $\R^n\subset(\Hat\C)^n$.  
Moreover, a right coset is $\Gamma(n)$-rational (resp. rational in the 
  $i$th coordinate) if and only if every coordinate (resp. the $i$th 
  coordinate) of the corresponding point in $\R^n$ is rational.  
It is totally irrational if and only if every coordinate of the 
  corresponding point in $\R^n$ is irrational.

\section{Encoding via piecewise linear functions}\label{S:PL}
In this section we associate a piecewise linear function $W_\bx:\R\to\R$ 
  to every $\bx\in\R^n$ and use it to characterise the set $E'_n$.  

We shall draw freely upon the standard results of continued fraction 
  theory.  All the results we use can be found in \cite{Kh}.  

For any $x\in\R$, let $W_x:\R\to\R$ be the piecewise linear function 
  determined by the conditions 
  \begin{enumerate}
    \item The function $W_x$ is continuous and nonnegative, 
    \item has slopes $\pm1$ whenever defined, 
    \item each local minimum of $W_x$ is a zero, and 
    \item the zeroes of $W_x$ are enumerated by $(2\log q_k)$ where 
     $(q_k)$ is the sequence of heights formed by the convergents of $x$.  
  \end{enumerate}
(Recall that the height of a rational is the smallest positive integer 
  that multiplies it into the integers.)  

It is well-known that an integer $q$ is the height of a convergent of $x$ 
  if and only if the distance of $qx$ to the nearest integer is minimal 
  among the set of integer multiples $x,2x,\dots,qx$.  
The function $W_x(t)$ is broken at a finite number of points if and only 
  if $x$ is rational.  

For any nonempty closed discrete subset $Z$ of $\R^2$ we define 
  $$\ell_\infty(Z):=\min\{\max(|a|,|b|): (a,b)\in Z, (a,b)\neq(0,0)\}$$ 
  to be the length of the shortest nonzero vector in $Z$ with respect 
  to the sup norm.  Let 
  $$g_t=\pmat{e^{t/2}&0\\0&e^{-t/2}}\quad\text{and}\quad h_s=\pmat{1&s\\0&1}.$$  
\begin{theorem}\label{thm:PL}
There is a universal constant $C>0$ such that 
     $$0\le|-2\log\ell_\infty(g_th_x^{-1}\Z^2)-W_x(t)|\le C$$ 
  for any $x\in\R$ and for any $t\in\R$.  
\end{theorem}
\begin{proof}
For any nonempty closed discrete subset $Z\subset\R^2$ the function 
  $$W(t)=-2\log\ell_\infty(g_tZ)$$ is a continuous piecewise linear 
  function with slopes $\pm1$ and isolated critical points, and in 
  particular, for the set $Z=h_x^{-1}\Z^2$.  
Observe that $W(t)$ has constant slope $-1$ for all $t<0$ since the 
  shortest nonzero vector in $g_tZ$ is realised by $(1,0)$.  

Suppose that $W(t)$ has a local minimum at a time $t\ge0$.  Then there 
  is a pair of vectors $v,v'\in Z$ and a square $S$ centered at the 
  origin containing $g_tv$ on one of its vertical edges and $g_tv'$ 
  on one of its horizontal edges such that $S$ contains no nonzero 
  vectors of $g_tZ$ in its interior.  
The vectors have the form $v=(qx-p,q)$ and $v'=(q'x-p',q')$ for some 
  $p,q,p',q'\in\Z$ and we may assume that $q\ge0$, $q'\ge1$, and if 
  $q'=1$ then $p'$ can be chosen so that $|q'x-p'|\le1/2$.  
The side of the square is $$2e^{-W(t)/2}=2e^{t/2}|qx-p|=2e^{-t/2}q'.$$  
By Minkowski's theorem, the area of the square is at most $4$ so that 
  $$e^{-W(t)}=q'|qx-p|\le1$$ from which it follows that $|qx-p|\le1/2$, 
  i.e. the distance to the nearest integer $||qx||=|qx-p|$.  
Since $g_{-t}S$ contains no nonzero vectors of $Z$ in its interior, 
  it follows that $q'$ is the height of some convergent of $x$.  
The corresponding zero of $W_x(t)$ occurs at time $2\log q'$ so that 
  $$W(t)-W_x(2\log q') = -\log q'||qx|| = t-2\log q'.$$  
Since $q\le q'$ and $|q'x-p'|\le|qx-p|$, the area of the parallelogram 
  spanned by $v$ and $v'$ is at most $2q'||qx||$; and since $Z$ is a 
  unimodular lattice, the area is at least one.  
Therefore, each local minimum of $W_x(t)$ is shifted upwards and to 
  the right relative to a zero of $W(t)$ by an equal amount of at 
  most $\log 2$ in both directions.  

Conversely, given any zero of $W(t)$, occuring at time $2\log q$, 
  there is a convergent $p/q$ of $x$ such that $(qx-p,q)\in Z$ lies 
  on the top of a rectangle $R$ symmetric with respect to the origin 
  and whose interior contains no nonzero vectors of $Z$.  
We may assume $R$ is chosen largest possible so that it contains a 
  point of $Z$ on one of its vertical edges.  
Since its horizontal are vertical sides are respectively at most one 
  and at least one, there is a (unique) time $t\ge0$ when $g_tR$ is 
  a square and $W(t)$ has a local minimum at this time.  

It follows easily that $$W_x(t)\le W(t)\le W_x(t)+2\log2$$ so that 
  the statement of the theorem holds with $C=2\log 2$.  
\end{proof}

\emph{Remark 1.}
The encoding of geodesics on the modular surface is well-known and dates 
  back to E.~Artin \cite{Ar}.  Theorem~\ref{thm:PL} may be thought of as 
  a refinement of this encoding that also involves the parametrisation 
  of the geodesics.  

\emph{Remark 2.} 
The set of oriented unimodular lattices in $\R^2$ may naturally be identified 
  with $G/\Gamma$ and the function $-2\log\ell_2(Z)$, where $\ell_2(Z)$ 
  denotes the Euclidean norm of the shortest nonzero vector in $Z$, lifts 
  to a $K$-invariant proper function on $G/\Gamma$.  The induced function 
  on the hyperbolic triangle $$\Delta:=\{z:\Im z>0, |\Re z|\le1/2, |z|\ge1\}$$ 
  coincides with the function $z\to \Im z$.  

For any $\bx\in\R^n$ we let $W_\bx(t):\R\to\R$ be the function given by 
  $$W_\bx(t) = \max (W_{x_1}(t),\dots,W_{x_n}(t)).$$  
Note that $W_\bx(t)$ has infinitely many local minima if and only if 
  every coordinate of $\bx$ is irrational.  
It follows by Mahler's criterion and Theorem~\ref{thm:PL} that 
  $\bx\in E_n$ if and only if $W_\bx(t)\to\infty$ as $t\to\infty$, 
  which in turn holds if and only if 
  $$W_\bx(t_j)\to\infty\quad\text{as}\quad j\to\infty$$ 
  where $(t_j)$ is the sequence of local minimum times.  
In particular, if $(q_k)$ and $(q'_l)$ are the sequences of heights 
  formed by the convergents of irrationals $x$ and $y$, respectively, 
  then $(x,y)\in E_2'$ if and only if 
  $$|\log(q_k/q'_l)|\to\infty\quad\text{as}\quad\min(q_k,q'_l)\to\infty;$$  
  in other words, for any $\rho>1$ there exists $R>1$ such that 
    for any indices $k,l$ such that $q_k>R$ and $q'_l>R$, we have 
    $$\max\left(\frac{q'_l}{q_k},\frac{q_k}{q'_l}\right)>\rho.$$

\section{Lower bound calculation}\label{S:Lower}
In this section we prove that $\Hdim E'_n\ge n-\frac{1}{2}$ for $n\ge2$.  
We begin by recalling some results that we need for this calculation.  

\textbf{Hausdorff dimension estimates.}
Consider a Cantor set $F\subset\R$ defined as a nested intersection 
  $$F=\bigcap_{j\ge0} F_j$$ where $F_0$ is a closed interval and each 
  $F_j$ is a disjoint union of (finitely) many closed intervals.  
Let $(m_j)_{j\ge1}$ be a sequence of positive integers and 
  $(\eps_j)_{j\ge1}$ a sequence of real numbers that tend to zero 
  monotonically and suppose that for each $j\ge1$ there are at least 
  $m_j$ intervals of $F_j$ contained in each interval of $F_{j-1}$ and 
  these intervals are separated by gaps of length at least $\eps_j$.  
Then the Hausdorff dimension of $F$ satisfies the lower bound 
  (see \cite{Fa}, Example~4.6) 
  \begin{equation}\label{ieq:lower}
    \Hdim F \ge \limsup_{j\to\infty} 
                 \frac{\log(m_1\cdots m_{j-1})}{-\log m_j\epsilon_j}.
  \end{equation}

The following product formula (see \cite{Fa}, Corollary~7.4) holds 
  for subsets $E\subset\R^n$ and $F\subset\R^m$,
  \begin{equation}\label{eqn:product} 
    \Hdim (E\times F) = \Hdim E + \Hdim F 
  \end{equation} 
  as soon as the Hausdorff and Minkowski dimensions of one of the 
  sets coincide.  
Thus, it is enough to prove $\Hdim E'_2 \ge \frac{3}{2}$.  

We shall also need the following result.  
\begin{lemma}\label{lem:projection}
(\cite{Fa}, Corollary~7.12)
Let $F\subset\R^2$ and $E$ its projection to the $x$-axis.  
For each $x\in E$, let $L_x$ be the line $\{(x,y):y\in\R\}$.  
If $\Hdim (F\cap L_x)\ge t$ for all $x\in E$ then 
  $$\Hdim F\ge t+\Hdim E.$$
\end{lemma}

\textbf{Counting rationals.}
We shall also need the following result.  
\begin{theorem}\label{thm:count}
(\cite{Ch1}, Theorem~2)
There is a constant $c>0$ such that for any interval $I$ of the form 
  $I=[x-d,x+d]$ and for any $h>0$ the number of rationals in $I$ whose 
  height lies between $h$ and $2h$ is at least $ch^2|I|$ provided $x$ 
  has a convergent whose height $q$ satisfies $(hd)^{-1}\le q\le h$.  
\end{theorem}
\textit{Remark.} Theorem~2 in \cite{Ch1} is stated under an addititonal 
  technical hypothesis that can be shown to be redundant and has therefore 
  been omitted here.  See \cite{Ch2} for other variations and improvements.  

We now show that $\Hdim E'_2\ge\frac{3}{2}$.  
The proof of the next lemma will be omitted, since it is essentially 
  contained in the proof of a classical result of Jarnik-Besicovitch.  
\begin{lemma}\label{lem:JB}
Let $A_\delta$ be the set of irrationals with the property that 
  the sequence $(q_k)$ of heights formed by its convergents satisfy 
  \begin{equation}\label{delta:speed}  
     q_k^{1+\delta} \le q_{k+1} \le 2q_k^{1+\delta}
  \end{equation}
  for $k$ large enough.  
Then $\Hdim A_\delta\ge\frac{1}{2+\delta}$ for any $\delta>0$.  
\end{lemma}

\begin{lemma}\label{lem:Hdim1}
Let $B_x$ be the set of $y\in\R$ such that $(x,y)\in E'_2$.  
Then $\Hdim B_x=1$ for any $x\in A_\delta$.  
\end{lemma}
\begin{proof}
Let $(q_k)$ be the sequence of heights formed by the convergents of 
  some given $x\in A_\delta$ and fix $k_0$ such that (\ref{delta:speed}) 
  holds for all $k\ge k_0$.  
We shall choose $k_0\gg1$ so that certain estimates that arise in the 
  course of the proof will hold.  
Our goal is to construct a set of $y\in B_x$ with the property that 
  for each $k\ge k_0$ there is a pair of consecutive convergents of $y$ 
  whose heights $q<q'$ satisfy (up to a constant factor) 
     $$\min\left(\frac{q'}{q_k},\frac{q_k}{q}\right)>\log q_k.$$
This set will be realised as a decreasing intersection $\cap F_j$ where 
  each $F_j$ is a disjoint union of closed intervals.  
Let us fix the following parameters for the construction: 
    $$ h_j  = \left[\frac{q_{k_0+j}}{\log q_{k_0+j}}\right],\quad 
     \eps_j = \frac{1}{8h_j^2}, \quad\text{and}\quad 
       d_j  = \frac{1}{q_{k_0+j}^2}$$ 
  where $[\cdot]$ denotes the greatest integer function.  
The intervals of $F_j$ will all have length $2d_j$ and the length of the 
  gap between any two of them will be at least $\eps_j$; moreover we shall 
  also require each interval of $F_j$ to be centered about a rational whose 
  height is between $h_j$ and $2h_j$.  

Let $F_0$ be an interval of length $2d_0$ centered about a rational of 
  height $h_0$.  
Given $F_j$ we shall define $F_{j+1}$ by specifying, for each interval 
  $I$ of $F_j$, the intervals of $F_{j+1}$ that are contained in $I$.  
Let $J$ be an interval of $F_j$ where $j\ge0$.  
By the induction hypothesis, $J$ is of the form 
      $$J=\left[\frac{p}{q}-d_j, \frac{p}{q}+d_j\right]$$ 
  where $p/q$ is a rational whose height $q$ is between $h_j$ and $2h_j$.  
Note that by construction, every $x\in J$ has $p/q$ as a convergent since 
  \begin{equation}\label{ieq:scc}
    \left|x-\frac{p}{q}\right| \le d_j = \frac{1}{q_{k_0+j}^2}
                               \le \frac{1}{8h_j^2} \le \frac{1}{2q^2}
  \end{equation}
  where in the third step we used $k_0\gg1$.  
Moreover, if $q'$ is the height of the next convergent then 
  \begin{equation}\label{ieq:cfi}
    \frac{1}{2qq'}<\left|x-\frac{p}{q}\right|<\frac{1}{qq'}.  
  \end{equation}
Let $I$ be the interval of length $\frac{d_j}{3}$ centered about 
  $$x'=\frac{p}{q} + \frac{d_j}{2}$$ and note that an interval of 
  length $2d_{j+1}$ centered about any point in $I$ is entirely contained 
  in $J$ provided $k_0\gg1$.  
The interval $I$ is of the form $[x'-d,x'+d]$ where $d=\frac{d_j}{6}$ 
  and the first inequality in (\ref{ieq:cfi}) implies 
  \begin{equation}\label{ieq:next}
    q' > \frac{1}{2qd_j} \ge \frac{1}{4h_jd_j} > (h_{j+1}d)^{-1} 
  \end{equation} 
  assuming $k_0\gg1$.  
For any point $y\in I$ the distance to $p/q$ is at least $\frac{d_j}{3}$ 
  so that the second inequality in (\ref{ieq:cfi}) implies 
  $$q' < \frac{3}{d_jq} \le \frac{3q_{k_0+j}^2}{h_j}
       < 6q_{k_0+j}\log q_{k_0+j} 
       < \frac{q_{k_0+j+1}}{\log q_{k_0+j+1}} \le h_{j+1}$$ 
  where in the fourth step we used $k_0\gg1$ again.  
Thus, Theorem~\ref{thm:count} applied to the interval $I$ implies there 
  are at least $$m_{j+1}=(c/3)h_{j+1}^2d_j$$ rationals in $I$ with heights 
  between $h_{j+1}$ and $2h_{j+1}$.  
We define the intervals of $F_{j+1}$ contained in $J$ to be any closed 
  interval of length $2d_{j+1}$ centered about a rational just constructed 
  and let $F=\cap F_j$ be the resulting Cantor set.  

Now we compute the Hausdorff dimension of $F$.  
The distance between the centers of two intervals of $F_j$ is at least 
  $2\eps_j$ because the height the rationals at the centers are both 
  at most $2h_j$.  
Since $2d_j<\eps_j$, the length of the gap between any two intervals 
  of $F_j$ is at least $\eps_j$.  
Before applying the formula (\ref{ieq:lower}), we need to develop some 
  estimates on the parameters of the construction.  
From (\ref{delta:speed}) we have 
  \begin{align*}
    q_{k_0}^{(1+\delta)^j} &\le q_{k_0+j}\le 2^jq_{k_0}^{(1+\delta)^j} \\
  \intertext{so that}
    \log q_{k_0+j} &= (1+\delta)^j\log q_{k_0} + O(j).  
  \end{align*}  
Let us write $A\asymp B$ if we have $A/C\le B\le AC$ for some explicitly 
  computable $C>0$ whose value is otherwise irrelevant for the purposes 
  of this calculation.  Then, using (\ref{delta:speed}) again, we have 
  \begin{align*}
    m_{j+1} &\asymp h_{j+1}^2d_j 
          \asymp \frac{q_{k_0+j}^{2\delta}}{(\log q_{k_0+j+1})^2} \\
  \intertext{so that}
    \log m_{j+1} &= 2\delta(1+\delta)^j\log q_{k_0} + O(j).  
  \end{align*}
Therefore, 
  \begin{align*}
    \log(m_1\cdots m_j) 
       &= \sum_{i=0}^{j-1}2\delta\log(1+\delta)^i\log q_{k_0} + O(j^2) \\
       &= 2(1+\delta)^j\log q_{k_0} + O(j^2)
  \end{align*}
  and since $m_{j+1}\eps_{j+1}\asymp d_j$ we have 
  $$-\log m_{j+1}\eps_{j+1}=2(1+\delta)^j\log q_{k_0} + O(j)$$ 
  so that the formula (\ref{ieq:lower}) yields $\Hdim F=1$.  

It remains to show that $F\subset B_x$.  
Given any $y\in F=\cap F_j$ and any $k=k_0+j$ for some $j\ge0$ we let 
  $p/q$ be the rational at the center of the interval $J$ of $F_j$ that 
  contains $y$.  Since $y\in J$, (\ref{ieq:scc}) implies $p/q$ is a 
  convergent of $y$ whose height $q$ satisfies 
              $$q<2h_j<\frac{2q_k}{\log q_k}.$$  
Now let $I\subset J$ be the interval of $F_{j+1}$ that contains $y$.  
Since $y\in I$, the first two inequalities in (\ref{ieq:next}) imply 
  $$q'>\frac{1}{4h_jd_j}>\frac{q_k\log q_k}{4}$$ 
  from which it follows that $(x,y)\in E_2'$.  
Therefore, $F\subset B_x$ and the lemma follows.  
\end{proof}

Lemmas~\ref{lem:projection}, \ref{lem:JB} and \ref{lem:Hdim1} now 
  imply that $\Hdim E'_2 \ge \frac{3}{2}$.

\section{Self-similar coverings}\label{S:Selfsim}
The main technical device for obtaining upper bounds on Hausdorff 
  dimension is the notion of a self-similar covering.  
\begin{definition}
Let $\cB$ be a countable covering of a subset $E\subset\R^n$ by 
  bounded subsets of $\R^n$ and let $\sigma$ be a function from 
  the set $\cB$ to the set of all nonempty subsets of $\cB$.  
We say $(\cB,\sigma)$ is a \emph{self-similar covering} of $E$ if 
  there exists $\lambda, 0<\lambda<1$ such that for every $x\in E$ 
  there is a sequence $B_j$ of elements in $\cB$ satisfying 
  \begin{enumerate}
    \item $\cap B_j=\{x\}$,
    \item $\diam B_{j+1}<\lambda \diam B_j$ for all $j$, and 
    \item $B_{j+1}\in\sigma(B_j)$ for all $j$.  
  \end{enumerate}
\end{definition}

We shall need a slightly more general notion which is notationally 
  more cumbersome but offers added flexibility in applications.  
\begin{definition}
Let $\cB$ be a countable covering of a subset $E\subset\R^n$ by 
  bounded subsets of $\R^n$ and assume that it is indexed by some 
  countable set $J$; let $\sigma$ be a function from the set $J$ 
  to the set of all nonempty subsets of $J$.  
For any $\alpha\in J$ we write $B(\alpha)$ for the element of $\cB$ 
  indexed by $\alpha$.  
We say $(\cB,J,\sigma)$ is an \emph{indexed self-similar covering} 
  of $E$ (the indexing function $\iota:J\to\cB$ being implicit) if 
  there exists a $\lambda, 0<\lambda<1$ such that for every $x\in E$ 
  we have a sequence $(\alpha_j)$ of elements in $J$ satisfying 
  \begin{enumerate}
    \item $\cap B(\alpha_j)=\{x\}$, 
    \item $\diam B(\alpha_{j+1})<\lambda \diam B(\alpha_j)$ 
            for all $j$, and 
    \item $\alpha_{j+1}\in\sigma(\alpha_j)$ for all $j$.  
  \end{enumerate}
\end{definition}

\textit{Remark.} 
We shall carry out the calculations using self-similar coverings by 
  indexed families of sets.  Although the calculations can also be done 
  using ordinary self-similar coverings, the advantage of using indexed 
  families will become clear, especially in the case $n>2$.  

The goal of this section is to prove the following.  
\begin{theorem}\label{thm:upper}
Let $(\cB,J,\sigma)$ be an indexed self-similar covering of a subset 
  $E\subset\R^n$ and suppose there is an $s>0$ such that for every 
  $\alpha\in J$ 
  \begin{equation}\label{ieq:selfsim}
    \sum_{\alpha'\in\sigma(\alpha)} 
                  (\diam B(\alpha'))^s\le (\diam B(\alpha))^s.
  \end{equation}  
Then $\Hdim E\le s$.  
\end{theorem}

For the convenience of the reader, we briefly recall the definition 
  of the Hausdorff dimension of a subset $E\subset\R^n$.  
A countable covering $\cU$ of a set $E\subset\R^n$ is said to be an 
  $\eps$-\emph{cover} if the diameter of each element is less than $\eps$.  
Given an $\eps$-cover $\cU$, we let $\mu^s(\cU)$ denote the sum of the 
  $s$th powers of the diameters of the elements in $\cU$, and set 
  $$\mu^s_\eps(E)=\inf
       \left\{\mu^s(\cU):\cU\text{ is an $\eps$-cover of $E$ }\right\}.$$  
The $s$-dimensional measure of $E$ is defined as the limit 
  \begin{equation*}
    \mu^s(E):=\lim_{\eps\to0^+}\mu^s_\eps(E) 
  \end{equation*}
  which as a function of $s\ge0$ has a critical value $h\ge0$ with the 
  property $\mu^s(E)=\infty$ for all $s<h$ while $\mu^s(E)=0$ for all $s>h$.  
This value may be taken as the definition of the Hausdorff dimension of $E$.  

The following is a basic property enjoyed by Hausdorff dimension.  
\begin{lemma}\label{lem:sup:Hdim}
(\cite{Bi}, Theorem~4.1)
For any countable collection $\{E_j\}$ of subsets of $\R^n$, 
              $$\Hdim (\cup E_j) = \sup \Hdim E_j.$$ 
\end{lemma}

\begin{proof}[Proof of Theorem~\ref{thm:upper}]
For each $x\in E$ we have a sequence of elements in $J$ satisfying 1.-3. 
  in the definition of a self-similar covering; we shall think of this 
  sequence as an infinite word $w(x)$ in a language over the alphabet $J$.  
If $w(x)=(\alpha_1,\alpha_2,\dots)$ then we say $(\alpha_1,\dots,\alpha_k)$ 
  is a \emph{prefix} of $w(x)$ of length $k$.  
Let $T$ be the collection of words formed by all possible prefixes of $w(x)$ 
  as $x$ ranges over $E$.  
For any $w\in T$, let $E_w$ be the subset of $E$ formed by those $x$ for 
  which the infinite word $w(x)$ has $w$ as a prefix.  
Since $T$ is countable and $E=\bigcup_{w\in T} E_w$ it suffices to show that 
  $\Hdim E_w\le s$ for every $w\in T$.  
For $w\in T$ we shall also write $B(w)$ for the element of $\cB$ indexed 
  by the last letter of $w$.  

Fix $w\in T$ and consider any $\eps>0$ with $\eps<\diam B(w)$.  
For each $x\in E_w$, let $j(x)$ be the smallest integer $j$ such that 
  $\diam B(\alpha_j(x))<\eps$ and let $w'(x)$ be the prefix of $w(x)$ 
  of length $j(x)$.   
Let $A\subset T$ be the set of all $w'(x)$ as $x$ ranges over $E_w$.  
Then $\cU=\{B(w'):w'\in A\}$ is an $\eps$-cover of $E_w$.  
We claim that $\mu^s(\cU)\le(\diam B(w))^s$ from which it follows that 
  $\mu^s_\eps(E_w)\le(\diam B(w))^s$ and since $\eps>0$ can be made 
  arbitrarily small, we have $\mu^s(E_w)<\infty$ so that $\Hdim E_w\le s$.  

To prove the claim, we consider the subset $T_w\subset T$ consisting 
  of all those words which have $w$ as a prefix.  
A subset $A'\subset T_w$ is an \emph{anti-chain} if no word in $A'$ 
  is a prefix of some other word in $A'$.  
By construction, $A$ is an anti-chain contained in $T_w$.  
Hence, the claim is a consequence of the following assertion: 
  \emph{for any anti-chain $A'\subset T_w$} 
  \begin{equation}\label{ieq:anti}
    \sum_{w'\in A'} (\diam B(w'))^s \le (\diam B(w))^s.  
  \end{equation}
Let $T_k\subset T_w$ be the subset formed by all words of length at 
  most $k+l$ where $l$ is the length of $w$.  
Since $T_0=\{w\}$, (\ref{ieq:anti}) holds for all $A'\subset T_0$.  
Now suppose that (\ref{ieq:anti}) holds for all anti-chains contained 
  in $T_k$ and let $A'$ be an anti-chain contained in $T_{k+1}$.  
Write $A'=A_0\cup A_1$ where $A_0=A\cap T_k$ and $A_1=A'- A_0$.  
Let $A''=A_0\cup\pi(A_1)$ where $\pi:A_1\to T_k$ is the map that sends 
  $w'$ to the word $\pi(w')$ obtained by dropping the last letter of $w'$.  
Note that $A_0\cap\pi(A_1)=\emptyset$ and that $A''$ is an anti-chain 
  contained in $T_k$.  If $w'\in\pi(A_1)$, $w''\in\pi^{-1}(w')$, and 
  $\alpha$ is the last letter of $w'$, then $w''$ is obtained from $w'$ 
  by adding a suffix $\alpha'\in\sigma(\alpha)$ where $\alpha$ is the 
  last letter of $w'$.  Moreover, it is obvious that $w''\in\pi(A_1)$ 
  is uniquely determined by $\alpha'$.  
Therefore, (using the notation $|\cdot|$ for the diameter of a set) 
\begin{align*}
  \sum_{w'\in A'} |B(w')|^s 
     &= \sum_{w'\in A_0} |B(w')|^s + 
          \sum_{w'\in\pi(A_1)} 
            \sum_{w''\in\pi^{-1}(w')} |B(w'')|^s \\
     &\le \sum_{w'\in A_0} |B(w')|^s + 
          \sum_{w'\in\pi(A_1)} 
            \sum_{\alpha'\in\sigma(\alpha)} |B(\alpha')|^s \\
     &\le \sum_{w'\in A''} |B(w')|^s 
\end{align*}
  which proves (\ref{ieq:anti}) for all anti-chains contained in $T_k$ 
  for some $k\ge0$.  
Since any $A'\subset T_w$ can be written as an increasing union of the 
  sets $A'_k=A'\cap T_k$, each of which is an anti-chain if $A'$ is, 
  it follows that (\ref{ieq:anti}) holds for all anti-chains contained 
  in $T_w$ as well, proving the assertion and hence the theorem.  
\end{proof}

As a special case of Theorem~\ref{thm:upper} we obtain.  
\begin{corollary}
Let $(\cB,\sigma)$ be a self-similar covering of a subset $E\subset\R^n$ 
  and suppose there is an $s>0$ such that for every $B\in\cB$ 
  \begin{equation*}
    \sum_{B'\in\sigma(B)} \left(\frac{\diam B'}{\diam B}\right)^s\le1.  
  \end{equation*}  
Then $\Hdim E\le s$.  
\end{corollary}

We remark that Theorem~\ref{thm:upper} is implied by the corollary in 
  the case when the indexing function $\iota:J\to\cB$ is injective.

\section{Upper bound calculation}\label{S:Upper}
In this section we prove the following.  
\begin{theorem}
Let $E_n(\delta),\delta>0$ be the set of points $\bx\in\R^n$ satisfying 
  \begin{enumerate}
    \item every coordinate of $\bx$ is irrational, and 
    \item $W_\bx(t)>-\log\delta$ for all sufficiently large $t$.  
  \end{enumerate}
Then for any $\delta\le2^{-7}$ we have 
  \begin{align*}
    \Hdim E_2(\delta) &\le \frac{3}{2}+16\delta + O(\delta^2)
  \intertext{and for $n>2$, there are constants $\delta_n>0$ and $C_n>0$ 
             such that} 
    \Hdim E_n(\delta) &\le n-\frac{1}{2}+C_n\sqrt{\delta} + O(\delta) 
  \end{align*}
  for any $\delta\le\delta_n$.  
\end{theorem}

Let $0<\delta<1$ be fixed.  
We shall first develop some convenient notation.  
Let $$Q=\{(p,q)\in\Z^2:\gcd(p,q)=1,q>0\}$$ 
  and for any $v,w\in Q$ such that $v=(p,q)$ and $w=(p',q')$ set 
  $$\Dot v=\frac{p}{q},\quad |v|=q,\quad v\times w=pq'-p'q.$$  

We shall introduce two relations $\vdash$ and $\models$ on $Q$ as follows.  
Let $\cC$ be the collection of sequences $\mathbf{v}:\N\cup\{0\}\to Q$ 
  with the property that $|v_0|=1$ and for each $k\ge0$, 
                   $$v_{k+1}=av_k+v_{k-1}$$ 
  for some $a\in\N$ with the understanding that when $k=0$ this holds 
  for some $a\ge2$ and some choice of $v_{-1}=\pm(1,0)$.  
Define $u\vdash v$ if there exists a sequence $(v_k)\in\cC$ and $k\ge0$ 
  such that $u=v_k$ and $v=v_{k+1}$; similarly, define $u\models v$ if 
  there exists a sequence $(v_k)\in\cC$, $k\ge0$ and $l\ge0$ such that 
  $u=v_k$ and $v=v_{k+l}$.  

It follows from well-known properties of continued fractions that a sequence 
  $(v_k)_{k\ge0}$ in $Q$ belongs to $\cC$ if and only if $(\Dot v_k)_{k\ge0}$ 
  is the sequence of convergents of some irrational number.  
It is easy to see that $|v_k|<|v_{k+1}|$ for all $k$ and for any sequence 
  $(v_k)$ in $\cC$; in other words the sequence of heights formed by the 
  convergents of any real number is strictly increasing.  
Note also that if $u\vdash v$ then $|u\times v|=1$ and $|u|<|v|$.  

\subsection{Case $n=2$}
For any $(u,v)\in Q\times Q$, let $B(u,v)$ denote the ball, with respect 
  to the sup metric on $\R^2$, of radius $\frac{1}{|u||v|}$ centered at 
  the rational point $(\Dot u,\Dot v)$.  
Let $\cB$ be the collection of balls $B(u,v)$ as $(u,v)$ ranges over 
  the elements of the set $J\subset Q\times Q$ consisting of all pairs 
  $(u,v)$ satisfying either of the mutually exclusive conditions 
          $$|u|<\delta|v| \quad\text{or}\quad |v|<\delta|u|.$$  
Given $(u,v)\in J$ satisfying $|u|<\delta|v|$ we define $\sigma(u,v)$ 
  to be the subset of $J$ consisting of all pairs $(u',v')$ satisfying 
  $$u\vdash u',\quad v\models v', \quad\text{and}\quad |v'|<\delta|u'|;$$  
  similarly, if $|v|<\delta|u|$ we define $\sigma(u,v)$ to be the subset 
  of $J$ consisting of those pairs $(u',v')$ satisfying 
  $v\vdash v'$, $u\models u'$, and $|u'|<\delta|v'|$.  

We now show $(\cB,J,\sigma)$ is a self-similar covering of $E_2(\delta)$.  
Given $\bx=(x,y)\in E_2(\delta)$ let $t_j\to\infty$ be the sequence of 
  times formed by the local minima of $W_\bx$ and choose any $j_0$ such 
  that $W_\bx(t_j)>-\log\delta$ for all $j>j_0$.  Let 
         $$\alpha_j=(u_j,v_j)\quad\text{for}\quad j>j_0$$  
  where $u_j,v_j\in Q$ are defined as follows.  
The zero of $W_x$ closest to $t_j$ is of the form $2\log q$ for some 
  integer $q$ which is the height of a convergent $p/q$ of $x$; let 
  $u_j$ be the element in $Q$ such that $\Dot u_j=p/q$.  
Similarly, let $v_j$ be the element of $Q$ such that $\Dot v_j$ is a 
  convergent of $y$ and $2\log|v_j|$ is the zero of $W_y$ closest to $t_j$.  

We claim that for any $j>j_0$   $$t_j=\log(|u_j||v_j|)\quad\text{and}\quad 
                 W_\bx(t_j)=\left|\log\frac{|u_j|}{|v_j|}\right|.$$ 
Since $W_\bx(t_j)>-\log\delta$ we either have $|u_j|<\delta|v_j|$ or 
  $|v_j|<\delta|u_j|$, so that the second part of the claim implies 
  $\alpha_j\in J$ for any $j>j_0$.  

To see the claim, let $t$ (resp. $t'$) be the time of the local maximum 
  of $W_\bx$ that immediately precedes (resp. follows) $t_j$, so that 
                         $$t<t_j<t'.$$  
Let $z\in\{x,y\}$ be the coordinate of $\bx$ such that $W_\bx(s)=W_z(s)$ 
  for all $s\in[t,t_j]$; similarly, let $z'$ be the coordinate of $\bx$ 
  such that $W_\bx(s)=W_{z'}(s)$ for all $s\in[t_j,t']$.  
Then $z$ has a consecutive pair of convergents with heights $p<p'$ such that 
     $$t=\log(pp')\quad\text{and}\quad W_\bx(t)=W_z(t)=\log\frac{p'}{p};$$ 
  similarly, $z'$ has a consecutive pair of convergents with heights $q<q'$ 
  such that 
     $$t'=\log(qq')\quad\text{and}\quad W_\bx(t')=W_{z'}(t)=\log\frac{q'}{q}.$$  
Since $W_z(t_j)=2\log p'-t_j$ and $W_{z'}(t_j)=t_j-2\log q$ and both are 
  equal to $W_\bx(t_j)$ we have 
     $$t_j=\log(p'q)\quad\text{and}\quad W_\bx(t_j)=\log\frac{p'}{q}.$$  
Now $W_z(t)>W_z(t_j)$ implies $p<q$ while $W_{z'}(t)>W_{z'}(t_j)$ implies 
  $q<q'$ and since $W_\bx(t_j)>-\log\delta>0$ we also have $q<p'$ so that 
  altogether we have the inequalities 
                         $$p<q<p'<q'.$$   
Since $t_j=\log(p'q)$ and there are no zeroes of $W_z$ lying strictly 
  between $2\log p$ and $2\log p'$, it follows that $2\log p'$ is the 
  zero of $W_z$ closest to $t_j$; similarly, $2\log q$ is the zero of 
  $W_{z'}$ closest of $t_j$.  
Moreover, since $p<q<p'<q'$ we cannot have $z=z'$ so that 
  $$\text{either}\quad \bx=(z,z')\quad\text{or}\quad \bx=(z',z).$$  
In the first case we have $|u_j|=p'$ and $|v_j|=q$, while in the second 
  case $|u_j|=q$ and $|v_j|=p'$; in both cases, the claim holds.  

Next, we verify that the sequence $(\alpha_j)_{j>j_0}$ satisfies the 
  three conditions in the definition of a self-similar covering.  
First, we show that $\bx\in B(\alpha_j)$ for all $j>j_0$.  
Suppose that $\alpha_j=(u_j,v_j)$ satisfies $|u_j|<\delta|v_j|$.  
Let $u'_j\in Q$ be the element such that $\Dot u'_j$ is the convergent 
  of $x$ that immediately follows $\Dot u_j$.  
Note that $$|v_j|<|u'_j|$$ since $t_j<t'$ where $t'=\log(|u_j||u'_j|)$ 
  is the time of the first local maximum of $W_\bx$ beyond $t_j$.  
It follows from continued fraction theory that 
  \begin{align*}
    |x-\Dot u_j| &< \frac{1}{|u_j||u_j'|}<\frac{1}{|u_j||v_j|}\\
  \intertext{and}
    |y-\Dot v_j| &< \frac{1}{|v_j|^2}<\frac{\delta}{|u_j||v_j|} 
  \end{align*}
  so that $\bx\in B(\alpha_j)$, assuming $|u_j|<\delta|v_j|$.  
If instead $\alpha_j=(u_j,v_j)$ satisfies $|v_j|<\delta|u_j|$, 
  an argument similar to the one just given leads to the same 
  conclusion $\bx\in B(\alpha_j)$.  
Therefore, $\bx\in\cap_{j>j_0}B(\alpha_j)$ and since 
  $\diam B(\alpha_j)=2e^{-2t_j}\to0$ as $j\to\infty$ it follows that 
  $\cap B(\alpha_j)=\{x\}$, giving the first condition in the definition 
  of a self-similar covering.  

Now we check the second and third conditions.  
Again, suppose that $\alpha_j=(u_j,v_j)$ satisfies $|u_j|<\delta|v_j|$.  
Let $z\in\{x,y\}$ the coordinate of $\bx$ such that $W_\bx(t)=W_z(t)$ 
  for all $t\in[t_j,t_{j+1}]$.  
Let $t'$ be the time of the unique local maximum of $W_\bx$ such that 
                        $$t_j<t'<t_{j+1}.$$ 
Then $t'=\log(qq')$ where $q<q'$ are the heights of a pair of consecutive 
  convergents of $z$.  From 
  \begin{align*}
    W_\bx(t')- W_\bx(t_j) &= t'-t_j,\\
               W_\bx(t_j) &= t_j-2\log|u_j|, \quad\text{and}\\
    W_\bx(t') = W_z(t')   &= t'-2\log q
  \end{align*}
  we see that $q=|u_j|$; and since $t_j<t'$ we also have $q'>|v_j|$.  
Thus $|v_j|$ lies strictly between $q$ and $q'$ so that $\Dot v_j$, 
  which is a convergent of $y$, cannot be a convergent of $z$; 
  it follows that $$z=x.$$  
And since $\Dot u_{j+1}$ is a convergent of $x$ with $|u_{j+1}|>|u_j|=q$ 
  and $q'$ is the height of the convergent that immediately follows 
  $\Dot u_j$, we have $q'\le|u_{j+1}|$.  
On the other hand, 
    $$2\log q'-t_{j+1}=W_\bx(t_{j+1})\ge\log\frac{|u_{j+1}|}{|v_{j+1}|}$$ 
  implies that $q'\ge|u_{j+1}|$, so that $q'=|u_{j+1}|$.  
It follows by the definition of $q$ and $q'$ that $u_j\vdash u_{j+1}$.  

We claim $|v_j|\le|v_{j+1}|$.  Indeed, suppose on the contrary that 
  $|v_{j+1}|<|v_j|$.  Let $q''$ be the height of the convergent of $y$ 
  that immediately follows $\Dot v_{j+1}$.  
Since $\Dot v_{j+1}$ precedes $\Dot v_j$ in the continued fraction 
  expansion of $y$, we must have $q''\le|v_j|$.  
Let $t''$ be the time of the local maximum of $W_\bx$ that immediately 
  follows $t_{j+1}$.  
Arguing as before, we see that $t''=\log(|v_{j+1}|q'')$ and since 
  $t''>t_{j+1}$ we have $q''>|u_{j+1}|=q'>|v_j|$, which gives a 
  contradiction, proving the claim.  

Since $v_j$ and $v_{j+1}$ are convergents of $y$, it follows that 
  $v_j\models v_{j+1}$, which together with $u_j\vdash u_{j+1}$ 
  establishes the third condition $$\alpha_{j+1}\in\sigma(\alpha_j)$$ 
  in the definition of a self-similar covering.  
From $$2\log q'-t_{j+1} = W(t_{j+1}) 
                        = \left|\log\frac{|u_{j+1}|}{|v_{j+1}|}\right|$$ 
  we see that $q'=\max(|u_{j+1}|,|v_{j+1}|)$.  
Since $q'=|u_{j+1}|$, as was shown earlier, we have $|v_{j+1}|\le|u_{j+1}|$ 
  and since $\alpha_{j+1}\in J$ we actually have $|v_{j+1}|<\delta|u_{j+1}|$.  
Therefore, $|u_j|<\delta|v_j|\le\delta|v_{j+1}|<\delta^2|u_{j+1}|$ so that 
  $$\frac{\diam B(\alpha_{j+1})}{\diam B(\alpha_j)}
                   =\frac{|u_j||v_j|}{|u_{j+1}||v_{j+1}|}<\delta^2<1$$ 
  giving the second condition in the definition of a self-similar covering 
  with $\lambda=\delta^2$.  
As the case where $\alpha_j=(u_j,v_j)$ satisfies $|v_j|<\delta|u_j|$ can be 
  treated by a similar argument, this completes the proof that 
  $(\cB,J,\sigma)$ is a self-similar covering of $E_2(\delta)$.  

Now we analyse the set $\sigma(u,v)$ for each $(u,v)\in J$.  
As before, we consider only the case $|u|<\delta|v|$ as the other case 
  is similar.  

Let $\pi:\sigma(u,v)\to\N\times\N$ be the function defined as follows.  
Given $(u',v')\in\sigma(u,v)$ we have $u\vdash u'$, $v\models v'$ and 
  $$|v'|<\delta|u'|.$$  
Since $u\vdash u'$ so that $|u\times u'|=1$ we have $$u'=au+\tu$$ for 
  some positive integer $a$ and some $\tu\in Q\cup\{\pm(1,0)\}$ satisfying 
  $|\tu|<|u|$; moreover, both $a$ and $\tu$ are uniquely determined by $u'$.  
Since $v\models v'$, there is a sequence $(v_k)\in\cC$ and $k,l\ge0$ such 
  that $v=v_k$ and $v'=v_{k+l}$.  We either have $v'=v$, $v'=v_{k+1}$, or 
  $v'$ is a positive linear combination of $v$ and $v_{k+1}$.  
Since $v_{k+1}=a'v+\tv$ for some uniquely determined integer $a'>0$ and 
  $\tv\in Q\cup\{\pm(1,0)\}$ such that $|\tv|<|v|$, we have in any case 
  $$v'=bv+c\tv$$ for some uniquely determined $b>0$ and $c\ge0$ such that 
  $c\le b$.  (In fact, $\gcd(b,c)=1$ although this will not be used.)  
We define $$\pi(u',v')=(a,b).$$  

For each $(u',v')\in\pi^{-1}(a,b)$ we have the inequalities 
  $$a|u|\le|u'|\le2a|u|\quad\text{and}\quad b|v|\le|v'|\le2b|v|$$ 
  so that 
  $$\frac{1}{4ab}\le\frac{\diam B(u',v')}{\diam B(u,v)}
                  =\frac{|u||v|}{|u'||v'|}\le\frac{1}{ab}.$$  
From the inequalities 
  $$|v|\le b|v|\le |v'|<\delta|u'|\le2a\delta|u|<2a\delta^2|v|$$ 
  we obtain the following bounds 
  \begin{equation*}
    a>\frac{1}{2\delta^2}\quad\text{and}\quad 
    1\le b\le 2a\delta^2.  
  \end{equation*}
Observe that $\pi^{-1}(a,b)$ contains at most $8b$ elements since 
  there are at most two choices for each of $\tu$ and $\tv$ while 
  there are at most $b+1\le2b$ choices for $c$.  

We now compute for any $s\in(1.5,2)$ 
\begin{align*}
  \sum_{(u',v')\in\sigma(u,v)}
      \left(\frac{\diam B(u',v')}{\diam B(u,v)}\right)^s
  &\le \sum_{a>1/2\delta^2}\sum_{b=1}^{[2a\delta^2]}\frac{8}{a^sb^{s-1}} \\
  &\le \frac{8\cdot(2\delta^2)^{2-s}}{2-s} 
       \sum_{a>1/2\delta^2}\frac{1}{a^{2s-2}} \\
  &\le \frac{8\cdot(2\delta^2)^{s-1}}{(2s-3)(2-s)} 
   \le \frac{16\delta}{(2s-3)(2-s)}.  
\end{align*}
The last expression is equal to one if $$2s^2-7s+6+16\delta=0$$ so that 
  Theorem~\ref{thm:upper} applies to give 
  $$\Hdim E_2(\delta)\le\frac{7}{4}-\sqrt{\frac{1}{16}-8\delta}
                       =\frac{3}{2}+16\delta+O(\delta^2)$$ 
  for any $\delta\le2^{-7}$.

\subsection{Case $n>2$}
Let $Q_n=Q^n\times\{1,\dots,n\}\times\{1,\dots,n\}$ and for any triple 
  $(\bc,i,j)\in Q_n$, let $B(\bc,i,j)$ denote the ball, with respect to 
  the sup metric on $\R^n$, of radius $\frac{1}{|c_i||c_j|}$ centered at 
  the rational point $\Dot\bc=(\Dot c_1,\dots,\Dot c_n)$.  
Let $\cB$ be the collection of balls $B(\bc,i,j)$ as $(\bc,i,j)$ ranges 
  over the elements of the set $J\subset Q_n$ consisting of all triples 
  $(\bc,i,j)$ satisfying 
  \begin{equation*}
    |c_j|<\sqrt{\delta}|c_i|\quad\text{and}\quad 
    \text{$|c_k|$ divides $|c_i||c_j|$ for $k=1,\dots,n$.}
  \end{equation*}  
For each $(\bc,i,j)\in Q_n$ let $A(\bc,i,j)\subset B(\bc,i,j)$ be the 
  subset consisting of those points whose $i$th coordinate lies within 
  a distance $\frac{1}{|c_i|^2}$ of $\Dot c_i$.  
For any $(\bc,i,j)\in\cB$ we define $\sigma(\bc,i,j)$ to be the subset 
  of $J$ consisting of all triples $(\bc',i',j')$ satisfying the following 
  conditions: 
  \begin{enumerate}
    \item[(a)] $i'=j$ and $c_j\vdash c'_j$, 
    \item[(b)] if $j'=i$ then $c_i\models c'_i$, 
    \item[(c)] $|c_i|<|c'_j|$ and $|c_j|<|c'_{j'}|$, 
    \item[(d)] $A(\bc,i,j)\cap A(\bc',i',j')\neq\emptyset$, and 
    \item[(e)] $|c_i|^2<|c'_{i'}||c'_{j'}|$.  
  \end{enumerate}

We now show $(\cB,J,\sigma)$ is a self-similar covering of $E_n(\delta)$.  
Given $\bx\in E_n(\delta)$ let $t_k\to\infty$ be the sequence of times 
  formed by the local maxima of $W_\bx$ and choose any $k_0$ such that 
  $W_\bx(t)>-\log\delta$ for all $t>t_{k_0}$.  
For each $k\ge k_0$ let 
             $$(u_k,v_k,i_k)\in Q\times Q\times\{1,\dots,n\}$$ 
  be a triple defined as follows.  
Let $i_k$ be the index $i$ such that $W_\bx(t_k)=W_{x_i}(t_k)$.  
Then $t_k=\log(qq')$ where $q<q'$ are heights of a pair of consecutive 
  convergents of $x_i$.  
Let $u_k$ (resp. $v_k$) be the element of $Q$ such that $\Dot u_k$ (resp. 
  $\Dot v_k$) is a convergent of $x_i$ and $|u_k|=q$ (resp. $|v_k|=q'$).  
For any $k>k_0$ we have 
              $$W_\bx(t)=\log\frac{|v_{k-1}|}{|u_k|}$$ 
  at the unique local minimum time $t$ between $t_{k-1}$ and $t_k$.  
Since $W_\bx(t)>-\log\delta$, we have $|u_k|<\delta|v_{k-1}|$ for all $k>k_0$.  

Construct a subsequence of $(u_k,v_k,i_k)$ as follows.  
Given $k_l$ let $k_{l+1}$ be the first (smallest) $k>k_l$ such that 
  \begin{equation}\label{def:k_l}
    |v_{k_l}|^2<|v_k||u_{k+1}| 
  \end{equation}
  holds.  
Let $(k_l)_{l\ge0}$ be the sequence obtained by recursive definition.  
By construction, we have 
  $$|v_{k_l}|^2<|v_{k_{l+1}}||u_{k_{l+2}}|$$ for all $l$.  
Moreover, since (\ref{def:k_l}) does not hold for $k=k_l$, we have 
  $$|v_{k_l}|^2\ge|v_{k_{l+1}-1}||u_{k_{l+1}}|>\delta^{-1}|u_{k_{l+1}}|^2$$ 
  for all $l$.  
To simply notation a bit, we shall suppress the double subscript and 
  write $(u_l,v_l,i_l)$ instead of $(u_{k_l},v_{k_l},i_{k_l})$.  
This will cause no confusion as we shall have no further use for the 
  original sequence of triples.  
Thus, for any $l\ge0$ the triple $(u_l,v_l,i_l)$ satisfies 
  \begin{equation}\label{def:ulvl}
    |u_{l+1}|<\sqrt{\delta}|v_l|,\qquad|v_l|^2<|v_{l+1}||u_{l+2}|,
  \end{equation} 
  and $\Dot u_l$ and $\Dot v_l$ are consecutive convergents of $x_{i_l}$.  

For any $l>0$ let $\bc_l$ be the element $(c_1,\dots,c_n)\in Q^n$ given by 
  \begin{equation*}
    c_i = \left\{\begin{tabular}{cc}
           $u_{l+1}$ & $i=i_{l+1}$ \\
           $v_l$     & $i=i_l$ \\
           $w_i$     & otherwise 
          \end{tabular}\right.
  \end{equation*}
  where $w_i\in Q$ is the element such that $\Dot w_i$ is the rational 
  closest to $x_i$ whose height $|w_i|$ divides $|u_{l+1}||v_l|$.  
Let $\alpha_l=(\bc_l,i_l,i_{l+1})$ and note that $\alpha_l\in J$ by the 
  first condition in (\ref{def:ulvl}) and the definition of $w_i$.  

Now we check that $(\alpha_l)_{l>0}$ satisfies the first condition in 
  the definition of a self-similar covering.  
To avoid double subscript notation, we fix $l>0$ and let 
  $(\bc,i,j)=\alpha_l=(\bc_l,i_l,i_{l+1})$.  Then 
  \begin{equation}\label{def:cicj}
     c_i=v_l,\qquad c_j=u_{l+1}
  \end{equation} 
  and for $k\neq i,j$ $\Dot c_k$ is the rational closest to $x_k$ with 
  height dividing $|c_i||c_j|$.  
Since $\Dot v_{l+1}$ is the convergent of $x_j$ that immediately 
  follows $\Dot c_j$, and $|v_l|<|v_{l+1}|$, we have 
  $$|x_j-\Dot c_j| < \frac{1}{|u_{l+1}||v_{l+1}|}<\frac{1}{|c_i||c_j|}.$$ 
Likewise, $\Dot c_i=\Dot v_l$ is a convergent of $x_i$ so that 
  $$|x_i-\Dot c_i| < \frac{1}{|v_l|^2}=\frac{1}{|c_i|^2}.$$  
By definition of $c_k$ we have (for $k\neq i,j$) 
  $$|x_k-\Dot c_k|\le\frac{1}{2|c_i||c_j|}<\frac{1}{|c_i||c_j|}.$$  
It follows that $\bx\in A(\alpha_l)$ for all $l>0$, and since 
  $A(\alpha_l)\subset B(\alpha_l)$ and 
  $\diam B(\alpha_l)=\frac{2}{|u_{l+1}||v_l|}\to0$ as $l\to\infty$ 
  we conclude that $\cap B(\alpha_l)=\{x\}$.  

Now we check that $(\alpha_l)$ satisfies the second and third conditions 
  in the definition of a self-similar covering.  
Fix $l>0$ and let $(\bc,i,j)=\alpha_l$ and $(\bc',i',j')=\alpha_{l+1}$.  
In addititon to (\ref{def:cicj}) we also have 
  $$i'=i_{l+l}=j,\quad c'_j=v_{l+1},\quad\text{and}\quad c'_{j'}=u_{l+2}.$$  
We first check that $\alpha_{l+1}$ satisfies the conditions (a)-(e) in 
  the definition of $\sigma(\alpha_l)$.  
Since $\Dot u_{l+1},\Dot v_{l+1}$ are consecutive convergents of $x_j$ and 
  $|u_{l+1}|<|v_{l+1}|$, we have $u_{l+1}\vdash v_{l+1}$ so that (a) follows.  
If $j'=i$ then $\Dot v_l,\Dot u_{l+2}$ are convergents of $x_i$ and since 
  $\Dot v_{l+2}$ is the convergent that immediately follows $\Dot u_{l+2}$, 
  we cannot have $|u_{l+2}|<|v_l|$ (because then $|v_{l+2}|\le|v_l|$, which 
  is a contradiction); therefore, $|v_l|\le|u_{l+2}|$ from which it follows 
  that $v_l\models u_{l+2}$, giving (b).  
The conditions in (c) follows from $|v_l|<|v_{l+1}|$ and $|u_{l+1}|<|u_{l+2}|$, 
  while (d) follows from $\bx\in A(\alpha_l)\cap A(\alpha_{l+1})$, and 
  (e) follows from the second condition in (\ref{def:ulvl}).  
Now using (e) and $\alpha_l\in J$ we find 
  $$\frac{\diam B(\alpha_{l+1})}{\diam B(\alpha_l)} 
          = \frac{|c_i||c_j|}{|c'_{i'}||c'_{j'}|}
          < \frac{|c_j|}{|c_i|} < \delta^{1/2} < 1$$
  so that the second condition in the definition of a self-similar covering 
  holds with $\lambda=\delta^{1/2}$.  

Now we analyse the set $\sigma(\alpha)$ for each $\alpha\in J$.  
Given $\alpha'\in\sigma(\alpha)$ we set $\alpha=(\bc,i,j)$ and 
  $\alpha'=(\bc',i',j')$ so that 
  $$\diam B(\alpha)=\frac{2}{|c_i||c_j|}\quad\text{and}\quad
    \diam B(\alpha')=\frac{2}{|c'_{i'}||c'_{j'}|}.$$  
Computing the sum in (\ref{ieq:selfsim}) reduces to analysing the 
  possibilities for the elements $c'_{i'}$ and $c'_{j'}$ of $Q$.  
The possibilities for $c'_{i'}$ can be analysed using the same techniques 
  in the case $n=2$; they can also be used to analyse the possibilities 
  for $c'_{j'}$ if $j'=i$.  
The main new ingredient is estimating the number of possibilities for 
  $c'_{j'}$ in the case when $j'\neq i$.  
Since $\diam B(\alpha')<\diam B(\alpha)$ we observe that $\Dot c'_{j'}$ 
  is a rational that lies within $\diam B(\alpha)$ of $\Dot c_{j'}$ and 
  is therefore restricted to some interval of length $2\diam B(\alpha)$ 
  centered at $c_{j'}$.  
Before continuing the calculation, we pause to state a lemma that will 
  be used to estimate the number of possibilities for $c'_{j'}$ in the 
  case when $j'\neq i$.  

\begin{lemma}\label{lem:count}
Let $I$ be an open interval of the real line and $b$ a positive integer.  
Suppose that the minimum height of a rational in $I$ is $q$.  
Then the number of rationals in $I$ with heights $q'$ satisfying 
  $\lfloor q'/q \rfloor=b$ is at most $2bq^2|I|$.  
\end{lemma}
\begin{proof}
A rational in $I$ satisfying the given condition corresponds 
to an integer point with relatively prime coordinates in the 
region $P:=\{(x,y)| x/y\in I, bq\leq y<(b+1)q\}$.  It suffices 
to show that $P$ contains at most $(b+1)q^2|I|$ integer points 
(relatively prime or not).  

Let $p/q$ be any rational in $I$ of minimum height.  
The projection of $P$ onto the $x$-axis along lines of 
slope $q/p$ is an open interval of length $(b+1)q|I|$.  
An integer point in $P$ maps to a rational whose height 
divides $q$.  Since a line of slope $q/p$ contains at 
most one integer point in $P$, it follows there are at 
most $(b+1)q^2|I|$ integer points in $P$.  
\end{proof}

Let $\sigma_k(\bc,i,j)$ be the subset of $\sigma(\bc,i,j)$ consisting 
  of those triples $(\bc',i',j')$ with $j'=k$.  
Note that $\sigma_j(\bc,i,j)=\emptyset$ since $j'\neq i'=j$.  
For $k\neq j$, let $\pi_k:\sigma_k(\bc,i,j)\to\N\times\N$ be the function 
  that sends $(\bc',j,k)$ to the pair $(a,b)$ given by 
  $$a = \left\lceil\frac{|c'_j|}{|c_j|}\right\rceil 
    \qquad\text{and}\qquad 
    b = \left\lfloor\frac{|c'_k|}{q_k}\right\rfloor$$  
  where $q_k:=|c_i|$ if $k=i$ and is otherwise set equal to the minimum 
  height of a rational in the open interval centered at $\Dot c_k$ of 
  length $\frac{4}{|c_i||c_j|}$.  
To see that $\pi_k$ is well-defined we need to check that $a\ge1$, which 
  holds by the second condition in (c), and that $b\ge1$, which, in the 
  case $k=i$, holds with equality and otherwise because the rational 
  $\Dot c'_k$ lies in the open interval centered at $\Dot c_k$ of length 
  $\frac{4}{|c_i||c_j|}$ and $q_k$ is the minimum height of a rational 
  in that interval.  

For each $(\bc',j,k)\in\pi_k^{-1}(a,b)$ we have the inequalities 
  $$|c'_j|\le a|c_j|\le 2|c'_j|\quad\text{and}\quad
   \frac{|c'_k|}{2}\le bq_k\le |c'_k|$$ 
  so that 
  \begin{equation}\label{ieq:ratio:n>2}
    \frac{|c_i|}{2abq_k} \le \frac{\diam B(\alpha')}{\diam B(\alpha)}
                           \le \frac{2|c_i|}{abq_k}.  
  \end{equation}
Moreover, the first condition in (c) and $\alpha\in J$ implies 
  $$a > \frac{|c'_j|}{|c_j|} > \frac{|c_i|}{|c_j|} > \frac{1}{\delta}$$ 
  while $\alpha'\in J$ and the first inequality above imply 
  $$|c'_k| < \sqrt{\delta}|c'_j| < a\sqrt{\delta}|c_j|$$ 
  giving us the following bounds 
  \begin{equation}\label{limits:n>2}
    a>\frac{1}{\sqrt{\delta}}\quad\text{and}\quad
    1\le b\le \frac{a\sqrt{\delta}|c_j|}{q_k}.  
  \end{equation}

Now we estimate the number of elements in $\pi_i^{-1}(a,b)$.  
Since each element is of the form $(\bc',j,i)$ we need to bound the 
  number of possibilities for $\bc'\in Q^n$.  
By (a) $c_j\vdash c'_j$ from which it follows that $c'_j=a'c_j+u$ 
  for some positive integer $a'$ and some $u\in Q\cup\{\pm(1,0)\}$ 
  satisfying $|c_j\times u|=1$ and $|u|<|c_j|$.  Note that $a'$ is 
  uniquely determined ($a'=a-1$) while $u$ admits two possibilities, 
  giving rise to $2$ choices for $c'_j$.  
By (b) $c_i\models c'_i$ from which it follows that $\Dot c'_i$ is 
  a convergent of $\Dot c_i$ and therefore is restricted to an open 
  interval of length $\frac{2}{|c_i|^2}$ centered at $\Dot c_i$.  
Since $\Dot c_i$ is the rational of minimum height in this interval 
  and $|c_i|=q_k$, Lemma~\ref{lem:count} implies there are at most 
  $4b$ choices for $c'_i$.  
For any index $k\neq i,j$ the rational $\Dot c'_k$ is restricted to 
  an open interval of length $\frac{4}{|c_i||c_j|}$ and since its 
  height $|c'_k|$ divides $|c'_i||c'_j|$, the number of choices for 
  $c'_k$ is bounded above by $\frac{4|c'_i||c'_j|}{|c_i||c_j|}$, which, 
  by the first inequality in (\ref{ieq:ratio:n>2}), is bounded by $8ab$.  
Therefore, $$^\#\pi_i^{-1}(a,b)\le8^{n-1}a^{n-2}b^{n-1}.$$  

Now we estimate the number of elements in $\pi_k^{-1}(a,b)$ for $k\neq i,j$.  
Again $i'=j$ and there are $2$ choices for $c'_j$.  
By (d) the rational $\Dot c'_k$ is restricted to an open interval of 
  length $\frac{4}{|c_i||c_j|}$ and since $q_k$ is the minimum height 
  of a rational in this interval, Lemma~\ref{lem:count} implies that 
  there are at most $\frac{8bq_k^2}{|c_i||c_j|}$ choices for $c'_k$.  
For any index $l\neq j,k$, the rational $\Dot c'_l$ is restricted to 
  an interval of length $\frac{4}{|c_i||c_j|}$ and since its height 
  $|c'_l|$ divides $|c'_j||c'_k|$, the number of choices for $c'_l$ 
  is bounded above by $\frac{8abq_k}{|c_i|}$.  
This estimate can be improved by a factor of $\frac{|c_j|}{|c_i|}$ in 
  the case $l=i$ because by (d) and (e) 
  $$|\Dot c'_i-\Dot c_i|<\frac{1}{|c'_{i'}||c'_{j'}|} + 
                         \frac{1}{|c_i|^2} \le \frac{2}{|c_i|^2}$$ 
  so that $\Dot c'_i$ is further restricted to an open interval whose 
  length is smaller than $2\diam B(\alpha)$ by this factor.  
Hence, for $k\neq i,j$, we have 
  $$^\#\pi_k^{-1}(a,b)\le 2\cdot8^{n-1}
             \left(\frac{q_k}{|c_i|}\right)^na^{n-2}b^{n-1}.$$  

Using (\ref{limits:n>2}) we now compute for any $k\neq i,j$, and 
  any $s\in(n-1/2,n)$ 
\begin{align*}
  \sum_{\alpha'\in\sigma_k(\alpha)}
      &\left(\frac{\diam B(\alpha')}{\diam B(\alpha)}\right)^s
   \le \sum_{a>\,\delta^{-1/2}}\sum_{b=1}^{a\sqrt{\delta}|c_j|q_k^{-1}}
     \sum_{\alpha'\in\pi_k^{-1}(a,b)}
      \left(\frac{2|c_i|}{abq_k}\right)^s \\
  &\le 2\cdot8^{n-1}\left(\frac{q_k}{|c_i|}\right)^{n-s} 
     \sum_{a>\,\delta^{-1/2}}\frac{1}{a^{s-n+2}}
     \sum_{b=1}^{a\sqrt{\delta}|c_j|q_k^{-1}}\frac{1}{b^{s-n+1}} \\
  &\le \frac{2\cdot8^{n-1}}{n-s}
       \left(\frac{\sqrt{\delta}|c_j|}{|c_i|}\right)^{n-s} 
     \sum_{a>\,\delta^{-1/2}}\frac{1}{a^{2s-2n+2}} \\
  &\le \frac{2\cdot8^{n-1}\sqrt{\delta}}{(2s-2n+1)(n-s)}
\end{align*}
  where we once again used $\alpha\in J$ in the last step.  
Note that if $k=i$ this estimate can be improved by a factor of $2$.  
Therefore, summing over $k\neq j$ we get 
\begin{align*}
  \sum_{\alpha'\in\sigma(\alpha)}
      &\left(\frac{\diam B(\alpha')}{\diam B(\alpha)}\right)^s
   \le \frac{(2n-3)\cdot8^{n-1}\sqrt{\delta}}{(2s-2n+1)(n-s)}.  
\end{align*}
Setting the last expression equal to one with $s=n-\frac{1}{2}+\eps$ 
  we find that $$2\eps^2-\eps+(2n-3)8^{n-1}\sqrt{\delta}=0$$ and 
  applying Theorem~\ref{thm:upper} we now obtain, for $n>2$, 
  $$\Hdim E_n(\delta) \le 
       n-\frac{1}{2} +(2n-3)\cdot8^{n-1}\sqrt{\delta} + O(\delta)$$
  for any $\delta\le 2^{-6n}(2n-3)^{-2}$.

\end{document}